\documentstyle[12pt]{article}
\makeatletter
\begin{document}
\thispagestyle{empty}
\pagestyle{empty}
\parskip = 10pt
\newcommand{\namelistlabel}[1]{\mbox{#1}\hfil}
\newenvironment{namelist}[1]{%
\begin{list}{}
{
\let\makelabel\namelistlabel
\settowidth{\labelwidth}{#1}
\setlength{\leftmargin}{1.1\labelwidth}
}
}{%
\end{list}}
\def\theequation{\thesection.\arabic{equation}}
\newtheorem{theorem}{\bf Theorem}[thetheorem]
\newtheorem{corollory}{\bf Corollory}[thecorollory]
\newtheorem{remark}{\bf Remark}[theremark]
\newtheorem{lemma}{\bf Lemma}[thelemma]
\def\thetheorem{\arabic{theorem}}
\def\theremark{\thesection.\arabic{remark}}
\def\thecorollory{\thesection.\arabic{corollory}}
\def\thelemma{\thesection.\arabic{lemma}}
\newcommand{\nc}{\newcommand}
\newcommand{\bsp}{\begin{sloppypar}}
\newcommand{\esp}{\end{sloppypar}}
\newcommand{\be}{\begin{equation}}
\newcommand{\ee}{\end{equation}}
\newcommand{\beanno}{\begin{eqnarray*}}
\newcommand{\inp}[2]{\left( {#1} ,\,{#2} \right)}
\newcommand{\dip}[2]{\left< {#1} ,\,{#2} \right>}
\newcommand{\disn}[1]{\|{#1}\|_h}
\newcommand{\pax}[1]{\frac{\partial{#1}}{\partial x}}
\newcommand{\tpar}[1]{\frac{\partial{#1}}{\partial t}}
\newcommand{\xpax}[2]{\frac{\partial^{#1}{#2}}{\partial x^{#1}}}
\newcommand{\pat}[2]{\frac{\partial^{#1}{#2}}{\partial t^{#1}}}
\newcommand{\ntpa}[2]{{\|\frac{\partial{#1}}{\partial t}\|}_{#2}}
\newcommand{\xpat}[2]{\frac{\partial^{#1}{#2}}{\partial t \partial x}}
\newcommand{\npat}[3]{{\|\frac{\partial^{#1}{#2}}{\partial t^{#1}}\|}_{#3}}
\newcommand{\eeanno}{\end{eqnarray*}}
\newcommand{\bea}{\begin{eqnarray}}
\newcommand{\eea}{\end{eqnarray}}
\newcommand{\ba}{\begin{array}}
\newcommand{\ea}{\end{array}}
\newcommand{\nno}{\nonumber}
\newcommand{\dou}{\partial}
\newcommand{\bc}{\begin{center}}
\newcommand{\ec}{\end{center}}
\newcommand{\bb}{\mbox{\hspace{.25cm}}}
\nc{\benu}{\begin{enumerate}}
\nc{\eenu}{\end{enumerate}}
\nc{\bth}{\begin{theorem}}
\nc{\eth}{\end{theorem}}
\nc{\bpr}{\begin{prop}}
\nc{\epr}{\end{prop}}
\nc{\blem}{\begin{lemma}}
\nc{\elem}{\end{lemma}}
\nc{\bcor}{\begin{corollary}}
\nc{\ecor}{\end{corollary}}
\newcommand{\R}{I\!\!\!R}
\nc{\Hy}{I\!\!\!H}
\nc{\C}{I\!\!\!C}
\newcommand{\la}{\lambda}
\newcommand{\s}{\sinh}
\newcommand{\co}{\cosh}
\newcommand{\vl}{V_{\la}}
\nc{\ga}{\gamma}
\nc{\vg}{V_{\ga}}
\nc{\T}{\Theta}
\nc{\f}{\frac}
\nc{\rw}{\rightarrow}
\nc{\om}{\omega}
\nc{\Om}{\Omega}
\nc{\al}{\alpha}
\nc{\qed}{\hfill \rule{2.5mm}{2.5mm}}
\setcounter{section}{0}
\title{Harmonic manifolds with some specific volume densities}                 
\author{ K.Ramachandran\thanks{Supported by National Board for Higher
Mathematics, DAE, INDIA}\,\,and A.Ranjan }
\maketitle
\begin{abstract}
We show that noncompact simply connected harmonic manifolds with volume
density $\Theta_{p}(r) =\sinh ^{n-1} r$ 
is isometric to the real hyperbolic space and noncompact simply connected
K\"{a}hler harmonic manifold 
with volume density $\Theta_{p}(r) =\sinh ^{2n-1} r \cosh r$ is isometric to 
the complex hyperbolic space. A similar result is also proved for
Quaternionic K\"{a}hler manifolds.
Using our methods we get an alternative
proof, without appealing to the powerful Cheeger-Gromoll splitting 
theorem, of the fact that every Ricci flat harmonic manifold is isometric
to the euclidean space. Finally a rigidity result for real hyperbolic
space is presented.
\end{abstract}
\section{Introduction}
\label{sec-intro}
\setcounter{equation}{0}
Let $(M,g)$ be a Riemannian manifold and let $p \in M$. Consider a normal 
coordinate neighbourhood $U$ around $p$. Let $\omega_{p}=\sqrt{|
\mbox{det}(g_{ij})|}$ be the volume density  function of $M$ in $U$.
We say that $M$ is a harmonic manifold if $\om_p$ is a function of
the geodesic distance $r(p,.)$ alone. If $(r,\phi)$ is a polar coordinate 
system around $p$ then the density becomes $\Theta_{p} = r^{n-1} \omega_{p}$.
So $M$ is harmonic if $\T_p$ is a function of $r$ alone and hence can be 
written as $\T_p(r)$. Moreover $\T_p(r)$ is independent of the point $p$
\cite{Bes}.

Rank one symmetric spaces are harmonic as can be easily seen from their
density function. Besides these there were no known examples of harmonic
spaces. Moreover Lichnerowicz proved that upto dimension 4 harmonic spaces
are in fact rank one symmetric. This led to the Lichnerowicz conjecture
which asserts that {\it Every harmonic space is rank one symmetric}. It
should be noted that even  higher rank symmetric spaces are not 
harmonic. 

Let $(M,g)$ be a harmonic space. The well known Ledger's formula
\cite{Bes} (see pp. 161)  gives
$$\om''_p(r)|_{r=0}\,=\,-\f{1}{3}\,Ricci_p$$
Hence for harmonic manifolds, since $\om_p$ is a function of $r$ alone,
the ricci curvature is a constant, i.e harmonic spaces are Einstein.
Let $Ricci(M) = k$. There arise three cases.\\
\begin{enumerate}
\item {\bf $k > 0.$}\, In this case, by Myers-Bonnet theorem, $M$ is 
compact with finite fundamental group and Szabo \cite{Sza} proved that 
compact harmonic  manifolds with  finite fundamental group 
are rank one symmetric, thus settling the Lichnerowicz conjecture.
\item {\bf $k = 0.$}\, Here one appeals to the powerful Cheeger-Gromoll
splitting theorem to conclude that $M$ is isometric to the euclidean
space, i.e {\it Ricci flat harmonic mmanifolds are flat.} 
\item {\bf $k < 0.$}\, The Lichnerowicz conjecture is not true in this case. 
E. Damek and F. Ricci \cite{DR} constructed a family of nonsymmetric
harmonic spaces. These spaces are called the $NA$ spaces. 
In this family there are harmonic manifolds with same 
density function as that of the quaternionic hyperbolic space. 
All these spaces are homogeneous. Presently it is not known whether
there are nonhomogeneous harmonic spaces. So it seems that the
classification of harmonic spaces can be achieved  only upto the 
determination of all density functions, i.e density  classification 
of harmonic spaces. So Szabo \cite{Sza1} asked the
following question.\\
{\bf Question.} Which harmonic spaces are determined by their density 
functions? i.e, which harmonic spaces are density equivalent?\\
\eenu
We answer this question for three specific cases, namely that of the real,
complex and Quaternionic hyperbolic spaces.

\begin{theorem}
Let $(M,g)$ be a non-compact simply connected harmonic space with density
function $\Theta_{p}(r) = \sinh^{n-1}r$, 
then $(M,g)$ is isometric to the real hyperbolic space. 
\end{theorem}
\begin{theorem}
Let (M,g) be a non-compact simply connected K\"{a}hler harmonic manifold
with density funtion $\Theta_{p}(r) =\sinh^{2n-1}r \cosh r $, 
then $M$ is isometric to the complex hyperbolic space.
\end{theorem}
\begin{theorem}
Let $(M,g)$ be a noncompact simply connected Quaternionic K\"{a}hler
harmonic manifold with volume density $\T_p(r) = \s^{4n-1}r \co^{3}r$,
then $M$ is isometric to the quaternionic hyperbolic space.
\end{theorem}
Theorems 1 and 2 explain the lack of examples of nonsymmetric (K\"{a}hler)
harmonic spaces with same density function as that of the real (complex)
hyperbolic space. Using the same methods 
we also give an alternative proof, without appealing to the 
powerful Cheeger-Gromoll splitting theorem, of the fact that every
ricci flat harmonic manifold is isometric to the euclidean space,i.e,
\begin{theorem}
Every simply connected ricci flat harmonic manifold is flat, i.e it is
isometric to the euclidean space.      
\end{theorem}

In the next section we give the proofs of the above theorems. In the last
section we give a rigidity result for the real hyperbolic space. 
The authors would like to thank their colleague  G. Santhanam for
discussions on the subject. 

\section{Density equivalent spaces}
\setcounter{equation}{0}
{\bf Proof of Theorem 1}\\ 
Let $p \in M$ and $S_{p,R}$ be the distance sphere around $p$ of radius $R$. 
From the well known Ledgers formula \cite{Bes} (see pp. 161) we have 
$$\nabla_{m} \nabla_{m} \omega_{p} = -2/3 Ricci(m,m) \;\; for \;\; p \in M
\;\;and \;\;m \in T_{p}(M).$$ 
hence we get $Ricci(g) = -(n-1)$. Again $ \Theta_p(r) = \sinh^{n-1} r$
gives that the mean curvature $\sigma_{p}(R)$ of $S_{p,R}$ is,
$$ \sigma_{p,R} = \Theta'_{p} / \Theta{p} = (n-1) \coth r.$$
Let $L$ be the second fundamental form of $S_{p,R}$, then $Trace L = (n-1)
\coth r$. Using the Riccati equation $ L'+ L^{2} + R\ (\gamma'\, .\, )\,
\gamma'\, = 0$ we get 
$$Tr L'+ TrL^{2} + Ricci = 0 \;\; i.e$$ 
$$Tr L^{2} = (n-1)\coth^{2} r = 1/(n-1)( Tr L)^{2}.$$
For any linear map $L$ we know that $ Tr L^{2} \geq 1/(n-1) (Tr L)^{2}$ and 
equality holds iff $L$ is a scalar operator. So $L$ is a scalar operator, i.e
$L = \coth r Id$ which shows that $M$ is of constant sectional curvature -1,
hence $M$ is isometric to the real hyperbolic space.
\qed

\vskip 1cm
\noindent
{\bf Proof of Theorem 2}\\
Ledgers formula gives $Ricci(g) = -(2n+2)$. 
Let $\gamma(t)$ be any geodesic. Let $J$ be the  complex structure on $M$ and
$I$ be the index form of $M$. Let $T>0$ be a real number. Let $ E_2 ,
E_3 , ... , E_{2n} $ be unit orthogonal parallel fields along $\gamma$,
normal to $\gamma'(t)$ with $ E_2 (t)\ =\ J\ \gamma'(t)$. Let $ J_i\ (t)$
be jacobi fields along $\gamma(t)$ such that
$$ J_i(0)\ =\ 0,\, \,  and\, \,  J_i(T)\ =\ E_i(T)$$
Let 
$$X_2(t)\,=\,\f{\s 2t}{\s 2T}\,E_2(t)$$
and
$$X_i(t)\ =\ \frac{\sinh t}{\sinh T} E_i(t),\,\,i\,=\,3,\cdots,2n. $$
be vector fields along $\gamma(t)$. Note that $ J_i(0) = X_i(0)$ and 
$J_i(T) = X_i(T)$. Since there are no conjugate points along $\gamma$, 
we get $ I(J_i,J_i)\ \leq\ I(X_i,X_i)$. Summing we get
$$ \sum _ {i=2}^{2n} I(J_i,J_i)\ \leq\ \sum _{i=2}^{2n}\ I(X_i,X_i)$$
But $ \sum_{i=2}^{2n} I(J_i,J_i)\,=\,\f{\Theta_{p}'(T)}{\Theta_{p}(T)}$.
A simple calculation gives
$$ I(X_2,X_2)\, =\, \frac{1}{2\sinh^2 2T}\, (\sinh 4T\, +\, 4T)\, -\, 
   \int_{0}^{T}\, \frac{\sinh^2 2t}{\sinh^2 2T}\, H(\gamma'(t))\,dt$$
$$ I(X_i,X_i)\,=\,  \frac{2T + \sinh 2T}{4 \sinh^2 T}\, -\,\int_0^T\, 
   \frac{\sinh^2\ t}{\sinh^2\ T}\ K(\gamma'(t) , E_i(t))\, dt\,\,,\,\,
i\,=\,3,\cdots,2n$$
Here $K(x,y)$ is the sectional curvature of the plane spanned by $x,y$,
and $H(\gamma'(t))\, =\, K(\gamma'(t),J\gamma'(t))$ is the holomorphic
sectional curvature.
Hence after simplifying and using Ricci = -2(n+1)
$$ \sum_{i=2}^{2n}\ I(X_i,X_i)\, =\ A(T)\ +\ \int_0^T\ B(t)\ H(\gamma'(t))\ dt $$
where 
$$ A(T)\, =\, \frac{\sinh 4T\, +\, 4T}{2\, \sinh^2\,2T}\, +\, 
   \frac{1}{4\,\sinh^2\,T}\, (4n\,\sinh 2T\, -\, 8T)$$
and
$$ B(t)\, =\,  \frac{\sinh^2\,t}{\sinh^2\, T}\, -\,
   \frac{\sinh^2\, 2t}{\sinh^2\, 2T}$$
Hence one gets 
$$ \frac{\Theta_p'(T)}{\Theta_p (T)}\, \leq\,( A(T)\, -\, 4\,C(T))\, +\, 
   \int_0^T\, B(t)\, (H(\gamma'(t))\, +\, 4)\, dt$$
where $C(T)\, =\, \int_0^T\, B(t)\,\,  dt$.  \\
Note that $B(t)\, \geq\, 0, t \in [0,T]$ and 
$A(T)\,-\,4\,C(T)\,=\,\frac{\Theta_p'\,(T)} {\Theta_p\,(T)}$. Hence it follows
that 
$$ \int_0^T\, B(t)\, (H(\gamma'(t))\,+\,4)\, dt\,\,\geq\,\,0$$
which in turn gives $H(v)\, \geq\, -4$ for all unit vectors $v$.
An algebraic calculation \cite{KN} yields 
$$ \int_{U_{p}M}\, H(v)\, dv\,=\,\frac{Vol(U_{p}M)}{n(n+1)}\, Scal_{p}M$$
where Scal is the scalar curvature of $M$. In our case it is $-4n(n+1)$.
So we get $$\int_{U_{p}M}\, H(v)\, dv\,=\,-4\, Vol(U_{p}M)$$
Combined with the conclusion $H(v) \geq -4$ we get $H\,\equiv\,-4$, and 
$M$ is isometric to the complex hyperbolic space.
\qed
\newpage
\noindent
{\bf Proof of Theorem 3}\\
Consider a chart $(U,p), p \in M$ with two almost complex $J_1,J_2$
such that the Levi-Civita derivatives of $J_1,J_2$ are linear
combinations of $J_1,J_2$ and $J_3 = J_1J_2$. Let $\ga(t)$ be a geodesic 
starting at $p$. Choose $T > 0$ such that $\ga[0,T] \subset U$. Let $I$ be
the index form on $\ga$. Let $E_2(t), \cdots, E_{4n}(t)$ be unit 
orthogonal parallel fields along $\ga$ such that $E_2(t), E_3(t), E_4(t)$
belong to the three dimensional subbundle spanned by $J_1 \ga'(t), 
J_2 \ga'(t)$ and $J_3 \ga'(t)$.

Now $Ricci(M) = -(4n+8)$ as can be seen from the Ledgers formula. 
The choice of $J_1, 
J_2, J_3$ shows that the four dimensional space spanned by 
$\{ \ga'(t), J_1 \ga'(t), J_2 \ga'(t), J_3 \ga'(t)$ is a quaternionic 
line parallel along $\ga$ which has an $Sp(1)$
action. Therefore we get a family of almost complex structures $J_1(t), 
J_2(t), J_3(t)$ along $\ga$ such that 
$$E_i(t)\,=\,J_i(t)\,,\,J_i(0)\,=\,J_i,\,{\rm
for}\,i\,=\,1,\,2,\,3$$ 
 
The rest of the proof is similar to that of theorem 2, so we shall be
brief. Take 
$$X_i(t)\,=\,\f{\s 2t}{\s 2T}\,E_i(t)\,,\,i\,=\,2,\,3,\,4$$
and
$$X_i(t)\,=\,\f{\s t}{\s T}\,E_i(t)\,,\,i\,=\,5,\,\cdots,\,4n$$
Let $J_i(t)$ be jacobi fields along $\ga(t)$ such that $J_i(0) = X_i(0)$
and $J_i(T) = X_i(T)$. Since there are no conjugate points along $\ga(t)$ 
the following inequality holds
$$ \sum _ {i=2}^{4n} I(J_i,J_i)\ \leq\ \sum _{i=2}^{4n}\ I(X_i,X_i)$$
and equality holds iff $X_i(t) = J_i(t)\,\,\forall\,i$. Now
$$\sum_{2}^{4n}\,I(X_i,X_i)\,=\,(A(T)\,-\,12\,C(T))\, +\, 
   \int_0^T\, B(t)\, \left(\sum_2^4\,K(\ga'(t),E_i(t))\, +\, 12\right)\, dt$$
where
$$A(T) = \f{3\,(\s 4T\,+\,4T)}{2\,\s^2 2T}\,+\,
\f{((8n\,+\,4)\,\s 2T\,-\,24T)}{4\,\s^2 T}$$
$$ B(t)\, =\,  \frac{\sinh^2\,t}{\sinh^2\, T}\, -\,
   \frac{\sinh^2\, 2t}{\sinh^2\, 2T}\,>\,0\,\,,\,t\,\in\,[0,T]$$
and
$$C(T)\, =\, \int_0^T\, B(t)\,\,  dt$$
and $K(x,y)$ stands for the sectional curvature of the plane spanned
by vectors $x$ and $y$. 
Using $\sum_{i=2}^{2n} I(J_i,J_i)\,=\,\f{\Theta_{p}'(T)}{\Theta_{p}(T)}$ 
and $E_i (t) = J_i (t)$ we get 
$$\f{\T'_p(T)}{\T_p(T)}\,\leq\,(A(T)\,-\,12\,C(T))\, +\, 
   \int_0^T\, B(t)\, \left(\sum_2^4\,K(\ga'(t),E_i(t))\, +\, 12\right)\, dt$$
Now we use the following relation between the components of the curvature
tensor \cite{Bes1},
$$\sum_1^3\,K\,(X\,,\,J_i(t)X)\,=\,\f{3}{n+2}\,Ricci(M)\,=\,-12$$
to finally get 
$$\f{\T'_p(T)}{\T_p(T)}\,\leq\,A(T)\,-\,12\,C(T)$$
and equality holds iff $X_i$ is a jacobi field for all $i$. A simple 
computation verifies that in fact equality holds. Thus $X_i(t) = J_i(t)$
for all $i$. Since the point $p$ is arbitrary and $\ga$ is any geodesic
starting at $p$, $M$ is isometric to the quaternionic
hyperbolic space. 
\qed
\vskip 1cm
\noindent
{\bf Proof of Theorem 4}\\
Let $\gamma$ be a geodesic ray in $M$. Let $F_{\gamma}$ be the Busemann
function relative to the geodesic ray $\gamma$.
Take $X_i(t)\, =\, \frac{t}{T}\, E_i(t)$ in the proof of Theorem 2 to get
$$\frac{\Theta_p'(T)}{\Theta_p(T)}\, \leq\, \frac{(n-1)}{T}$$
But convexity of balls gives $\f{\T'_p(T)}{\T_p(T)} \geq 0$. 
Therefore 
$$ \frac{\Theta_p'(T)}{\Theta_p(T)}\,\,\rightarrow\, 0\,\, as\,\,
T\,\rightarrow\, \infty$$ 
i.e, $Tr L\,=\,0$ 
where $L$ is the {\em second fundamental form} of the horospheres
determined by $F_{\gamma}$. Since Ricci = 0, Riccati equation gives
$Tr\,L'\,+\,Tr\,(L^2)\,=\,0$ , but $Tr\,L'\,=\,0$, hence 
$Tr\,(L^2)\,=\,0$. Symmetry of $L$ now gives that $L\,=\,0$. This shows 
that the sectional curvature $K(\ga'(t),.) = 0$. But $\ga$ is a arbitrary
geodesic, hence $K(M) \equiv 0$  and the proof is complete. 
\qed

\section{Rigidity of $\Hy^n$}
In this section we prove that the real hyperbolic space is rigid among all
harmonic spaces. More precisely we show that, if a harmonic space is
asymptotically density equivalent to the real hyperbolic space then 
they are actually isometric.

Let $(M,g)$ be a non-compact harmonic manifold.
Normalising the metric on $(M,g)$, let us assume that
$Ricci(M,g)\,=\,-(n-1)$. The Bishop-Gromov volume comparison theorem
\cite{GHL} (pp. 144-147) or \cite{Kum} (pp. 140)  gives the density
function to be 
$$\Theta_p(r)\,=\,\alpha(r)\,\sinh^{n-1} r,$$
where $ \alpha(r)$  satisfies $$\,0\,\leq\,\alpha(r)\,\leq\,1\,;\,
\alpha(0)\,=\,1\,,\,\alpha'(0)\,=\,0\,{\rm and}\,\alpha'(r)\,\leq\,0.$$ 
Hence $\alpha(r)$ is a decreasing function. Two cases arise:\\
1. $\lim_{r\,\rightarrow\,\infty}\,\alpha(r)\,=\,0$ and\\
2. $\lim_{r\,\rightarrow\,\infty}\,\alpha(r)\,=\,c$ for some constant $c >
0.$\\
In the second case the density of $M$ is asymptotically same as that of
the real hyperbolic case. In this case we show that $c = 1$ and $M$ is
in fact isometric to the real hyperbolic space.
\begin{theorem}
Let $(M,g)$ be a non-compact harmonic manifold with
$\lim_{r\,\rightarrow\,\infty}\,\alpha(r)\,=\,c\,>\,0 .$ Then $c = 1$ and
$M$ is isometric to the real hyperbolic space of constant sectional
curvature -1. 
\end{theorem}
{\bf Proof.} 
Let $\gamma$ be a geodesic ray in $M$. Let $F_{\gamma}$ be the Busemann
function of $\gamma$. Now
 $$\lim_{r\,\rightarrow\,\infty}\,\alpha(r)\,=\,c\,>\,0$$
implies that the density of the horospheres determined by $F_{\gamma}$ is
$$\Theta\,(r)\, =\, c\,\sinh^{n-1} r$$
Hence the mean curvature of these horospheres is
            $$\sigma_p\,(r)\,=\,(n-1)\,\coth r\,=\,Tr\,L$$
where $L$ is the {\em second fundamental form} of the horospheres. Now the
Riccati equation combined with  $Ricci\,=\,-(n-1)$ gives
$$Tr\,(L^2)\,=\,(n-1)\,\coth^2 r$$ but $$Tr\,L\,=\,(n-1)\,\coth r$$
Hence by Cauchy-Schwartz inequality one gets that $L\,=\,(\coth r)\,Id$.
Thus $M$ is isometric to the real hyperbolic space of constant curvature
-1. 
\qed
\section{Remarks}
For the complex hyperbolic space an easy calculation shows that
$\al(r) \rw 0$ as $r \rw \infty$. Now assume that $M$ is K\"ahler
harmonic space. Normalize the metric so that $Ricci(M) = -(2n+2)$.
Again applying the Bishop-Gromov comparison theorem one sees that
the volume density of $M$ is 
$$\T(r)\,=\,\beta(r)\,\s^{2n-1} r\,\co r$$
where $\beta(r)$ satisfies the same properties as that of $\al(r)$.
The following question is natural.\\
{\bf Question.} If $\beta(r) \rw b(>0)$ as $r \rw \infty$, is 
$M$ isometric to the complex hyperbolic space. 

\bibliographystyle{plain}

\begin{thebibliography}{99}
\bibitem{Bes} A. L. Besse, {\em Manifolds all of whose geodesics are
closed}, Springer, Berlin, 1978. 
\bibitem{Bes1} A. L. Besse, {\em Einstein Manifolds}, Springer, Berlin, 
1987.
\bibitem{DR} E. Damek and F. Ricci, {\em A class of nonsymmetric harmonic
riemmanian spaces}, Bulletin AMS, 27(1992), pp. 139-142.
\bibitem{GHL} S. Gallot, D. Hulin, J. Lafontaine, {\em Riemannian
Geometry}, Springer, 1990.
\bibitem{KN} S. Kobayashi, K. Nomizu, {\em Foundations of Differential 
Geometry, I, II}, Interscience, Wiley, Newyork, 1963, 1969.
\bibitem{Kum} S. Kumaresan, {\em A Course in Riemannian Geometry},
to appear.
\bibitem{Sza} Z.I.Szabo, {\em The Lichnerowicz conjecture on Harmonic
Manifolds}, J. Differential Geometry, 31(1990), pp. 1-28. 
\bibitem{Sza1} Z.I.Szabo, {\em Spectral theory for operator families on
Riemannian Manifolds}, Proceedings of symposia in Pure Math, 54 (1993)
part 3, pp. 615-665.
\end{thebibliography}

\noindent
Department of Mathematics,    \hspace*{3cm}kram@ganit.math.iitb.ernet.in\\
I.I.T Powai, Mumbai 400 076, \hspace*{2.85cm}aranjan@ganit.math.iitb.ernet.in\\
India.\\
\end{document}